\documentclass[12pt]{amsart}
\usepackage{amscd}      % to be able to use "CD" environment
\usepackage{amssymb}
\usepackage{amsmath, amsthm, graphics}
\usepackage{xypic}      % to be able to use "diagram" environment of
\LaTeXdiagrams          % xypic to typeset commutative diagrams
\usepackage[all,v2]{xy}
\xyoption{2cell} \UseAllTwocells \xyoption{frame} \CompileMatrices
\allowdisplaybreaks[3]

1

\usepackage{latexsym}
\usepackage{epsfig}
\usepackage{amsfonts}
\usepackage{enumerate}
\usepackage{times}

\theoremstyle{definition}

\theoremstyle{remark}

\theoremstyle{remark}

\def\<{\left\langle}
\def\>{\right\rangle}
\begin{document}

\title{Descendent bounds for effective divisors on $\overline{M}_g$}

\begin{abstract} The slope 
of $\overline{M}_g$ is bounded from below by $\frac{60}{g+4}$
via a  descendent calculation.
\end{abstract}

\author{R. Pandharipande}\thanks{Partially supported by DMS-0500187}

\date{\today}

\maketitle

\section{Slope}
Let $\overline{M}_g$ be the moduli space of curves for $g\geq 2$. Let
$\lambda \in A^1(\overline{M}_g)$
be the first Chern class of the Hodge bundle. Let
$\delta_0\in   A^1(\overline{M}_g)$ be the class of the
boundary divisor of irreducible
nodal curves. For 
$1 \leq i \leq \lfloor \frac{g}{2} \rfloor$,
let $\delta_i\in   A^1(\overline{M}_g)$ be the class of the
corresponding reducible boundary divisor.
We will consider {\em effective} divisors $D\subset \overline{M}_g$
of the form
\begin{equation} \label{ccc}
[D] = 
\alpha \lambda - \beta_0 \delta_0 - \sum_{i=1}^{\lfloor \frac{g}{2} \rfloor}
\beta_i \delta_i\ , \ \ 
\  \ \ \ \alpha,\beta_i \in 
\mathbb{Q}_{> 0}.
\end{equation}
Our goal is to find lower bounds for the slope $\frac{\alpha}{\beta_0}$
of $D$.

If $C\subset \overline{M}_g$ is a curve with no components contained in
$D$, then  
$$[C]\cdot [D] \geq 0.$$
If also $[C]\cdot \lambda >0$ and $C$ has no
components contained in the boundary divisors, 
then we conclude
$$\frac{\alpha}{\beta_0} \geq \frac{[C]\cdot \delta_0}{[C]\cdot
\lambda}.$$
Hence, such curves $C$ provide lower bounds for the slope.

\section{Cotangent lines}
Let $\psi_i \in A^1(\overline{M}_{g,n})$  be the first
Chern class of the 
cotangent line bundle
$$L_i \rightarrow \overline{M}_{g,n}$$
associated to the $i^{th}$ marked point.
The line bundles $L_i$ 
are well-known to be nef via the Hodge index
theorem for surfaces. Hence the 
class
\begin{equation}\label{vvv}
\psi_1^{r_1} \cdots \psi_n^{r_n} \in A_1(\overline{M}_{g,n}), \ \ \
\sum_{i=1}^n r_i = 3g-3+n-1
\end{equation}
is a limit of effective curve classes which sweep $\overline{M}_{g,n}$.
We can push the curve class \eqref{vvv} forward by
$$\epsilon: \overline{M}_{g,n} \rightarrow \overline{M}_g.$$
Therefore,
\begin{equation}
\label{nnn}
\frac{ \int_{\overline{M}_{g,n}} \psi_1^{r_1} \cdots \psi_n^{r_n} 
\cdot \delta_0}{
\int_{\overline{M}_{g,n}} \psi_1^{r_1} \cdots \psi_n^{r_n} 
\cdot \lambda}
\end{equation}
is a lower bound for slopes of effective divisors $D$ of the
form \eqref{ccc} as long as the denominator is positive.

\section{Calculation}
We calculate the bound  \eqref{nnn} in case $n=1$ and $r_1=3g-3$ by 
explicit evaluation of the two integrals.

Consider first the numerator. Using the normalization map of the
irreducible divisor,
$$\int_{\overline{M}_{g,1}} \psi_1^{3g-3} 
\cdot \delta_0 = \frac{1}{2}\int_{\overline{M}_{g-1,3}} 
\psi_1^{3g-3} =
\frac{1}{2}
\int_{\overline{M}_{g-1,1}} 
\psi_1^{3g-5}.$$
The string equation is used in the last equality.
The evaluation
$
\int_{\overline{M}_{g-1,1}} 
\psi_1^{3g-5} = \frac{1}{(24)^{g-1}(g-1)!}$ is well-known \cite{fp2}. Hence
$$\int_{\overline{M}_{g,1}} \psi_1^{3g-3} 
\cdot \delta_0 =\frac{1}{2}\frac{1}{(24)^{g-1}(g-1)!}.$$

The denominator is more complicated to evaluate. The first step
is to use the GRR equation of \cite{fp} for Hodge integrals,
\begin{equation}\label{dqq}
\int_{\overline{M}_{g,1}} \psi_1^{3g-3} 
\cdot \lambda
= \frac{B_2}{2} \int_{\overline{M}_{g,2}} \psi_1^{3g-3} \psi_2^2
-\frac{B_2}{2} \int_{\overline{M}_{g,1}} \psi_1^{3g-2} 
+ \frac{B_2}{4} \int_{\overline{M}_{g-1,3}} \psi_1^{3g-3}, 
\end{equation}
where $B_2=1/6$ is the second Bernoulli number.
We have already seen how to evaluate the last two integrals.
The first integral on the right side of \eqref{dqq} is evaluated
using the $L_1$ Virasoro constraint \cite{Wit},
$$\frac{15}{4}\int_{\overline{M}_{g,2}} \psi_1^{3g-3} \psi_2^2
= \frac{(6g-5)(6g-3)}{4} \int_{\overline{M}_{g,1}} \psi_1^{3g-2}
+ \frac{1}{2}\cdot \frac{1}{4}\int_{\overline{M}_{g-1,3}} \psi_1^{3g-3}.
$$
Putting this together, we find
$$\int_{\overline{M}_{g,1}} \psi_1^{3g-3} 
\cdot \lambda
= \left(\frac{(6g-5)(6g-3)}{180}-\frac{1}{12} \right) \frac{1}{(24)^{g}g!}
+\frac{2}{45} \frac{1}{(24)^{g-1}(g-1)!}.$$

After taking the ratio, we obtain a simple exact evaluation of
 the bound,
$$\frac{\int_{\overline{M}_{g,1}} \psi_1^{3g-3} 
\cdot \delta_0}{
\int_{\overline{M}_{g,1}} \psi_1^{3g-3} 
\cdot \lambda} = \frac{60}{g+4}.$$

\vspace{+10pt}

\noindent{\bf Proposition.} {\em The slope of $D\subset \overline{M}_g$
is always at least $\frac{60}{g+4}$.}

\vspace{+10pt}

Taking the limit for large $g$, we find
$$\lim_{g\rightarrow \infty}\left(
\frac{ \int_{\overline{M}_{g,1}} \psi_1^{3g-3}  
\cdot \delta_0}{
\int_{\overline{M}_{g,1}} \psi_1^{3g-3} 
\cdot \lambda}\right) \sim 
 \frac{60}{g}\ .$$
Hence, we derive the asymptotic $\frac{1}{g}$ bound
 predicted experimentally in \cite{hm}.

The method of \cite{hm} is the same, but the moving curves
there are obtained from Hurwitz covers rather than 
descendent intersections. The combinatorics of Hurwitz covers makes
exact analysis difficult, but computer calculations \cite{hm} 
predict the slope is always at least $\frac{576}{5g}$,
which is not very different from the Proposition.
The advantage of the descendent approach is the calculational
simplicity.

\section{Other bounds}
Low genus computations via Faber's program suggest the following property always holds.

\noindent{\bf Conjecture.}
For $g\geq 1$ and $\sum_{i=1}^n r_i = 3g-3+n-1$
$$\frac{\int_{\overline{M}_{g,1}} \psi_1^{3g-3} 
\cdot \delta_0}{
\int_{\overline{M}_{g,1}} \psi_1^{3g-3} 
\cdot \lambda} \geq 
\frac{ \int_{\overline{M}_{g,n}} \psi_1^{r_1} \cdots \psi_n^{r_n} 
\cdot \delta_0}{
\int_{\overline{M}_{g,n}} \psi_1^{r_1} \cdots \psi_n^{r_n} 
\cdot \lambda}.$$

If true, the bound of the Proposition is the best obtainable
from descendent integrals on the moduli space of curves.

\section{Acknowledgments}
The calculation was done for D. Chen who asked if any exact
bounds for slopes were possible to obtain. G. Farkas
has found a related approach using Weil-Petersson volumes.

\end{document}